

\input amstex
\documentstyle{amsppt}

\magnification=1200
\parskip 12pt
\pagewidth{5.4in}
\pageheight{7.2in}

\baselineskip=12pt
\expandafter\redefine\csname logo\string@\endcsname{}
\NoBlackBoxes                
\NoRunningHeads
\redefine\no{\noindent}

\define\C{\bold C}

\define\Z{\bold Z}

\redefine\P{\bold P}

\define\al{\alpha}

\define\de{\delta}
\define\la{\lambda}
 
\define\Si{\Sigma}

\define\om{\omega}

\define\sub{\subseteq}

\define\st{\ \vert\ }   

\redefine\ll{\lq\lq}
\redefine\rr{\rq\rq\ }
\define\rrr{\rq\rq}

\redefine\dim{\operatorname {dim}}

\define\End{\operatorname {End}}

\define\lan{\langle}
\define\ran{\rangle}
\define\llan{\lan\lan}
\define\rran{\ran\ran}

\define\subD{_{\ssize D}}
\define\HolD{\operatorname {Hol}\subD}


\define\ct{\circ_{t}}
\define\cth{\circ_{\hat t}}

\redefine\b{\partial}
\define\bone{\b_1}
\define\btwo{\b_2}
\define\bi{\b_i}
\redefine\br{\b_r}

\define\ddtone{\frac{d}{dt_1}}
\define\ddqone{\frac{d}{dq_1}}

\define\ddone{\frac{\b}{\b t_1}}
\define\ddare{\frac{\b}{\b t_r}}
\define\ddi{\frac{\b}{\b t_i}}
\define\ddj{\frac{\b}{\b t_j}}

\define\plus{\ +\ }
\redefine\llan{\lan\hskip-2pt\lan}
\redefine\rran{\ran\hskip-2pt\ran}

\define\lkr{{\ \leftarrow k \rightarrow\ }}
\define\lkmor{{\ \leftarrow k-1 \rightarrow\ }}

\topmatter
\title Introduction to homological geometry: part II
\endtitle
\author Martin A. Guest
\endauthor
\thanks The author was partially supported by a Grant-in-Aid for Scientific Research
from the Ministry of Education.
\endthanks
\endtopmatter

\document

Quantum cohomology is a concrete
manifestation of the deep relations between topology and integrable
systems which have been suggested by quantum field theory.  It is a generalization of
ordinary cohomology theory, and (thanks to the pioneering work of
many people) it can be defined in a purely mathematical
way, but it has two striking features.  First, there is as
yet no general set of computational techniques analogous to the standard machinery of
algebraic topology.  Many calculations have been carried out for individual
spaces, but few general methods of calculation are known.  Second,
because of the link with integrable systems, a fundamental role is
played by various functions (differential
operators, connections, etc.), and the variables involved in these functions
are homology or cohomology classes.  It seems likely that this
function-theoretic point of view will form the basis of the desired
machinery. In fact, major progress in this direction
has been initiated by A. Givental, who has invented a new term for it,
\ll homological geometry\rrr.  

It is easy to explain the origin of these homological functions:  a cohomology
theory amounts to a discrete collection of data which expresses the possible intersections
of the various homology classes for a given manifold $M$, and the functions
in question serve as generating functions for this morass of combinatorial information.
For ordinary cohomology, the data is finite, but for quantum cohomology
it is usually not.  Therefore the problem is one of dealing expeditiously with
this collection of numbers, which includes, in particular,
the \ll structure constants\rr
of the quantum cohomology ring.  These numbers are called Gromov-Witten invariants.

As the authors of the textbook \cite{GKP} state
(at the beginning of chapter 7), \ll The most powerful way to deal
with sequences of numbers, as far as anybody knows, is to manipulate
infinite series that generate those sequences\rrr. Accordingly, generating
functions form the basis of the function-theoretic approach to quantum cohomology. It is
particularly gratifying when these generating functions, originally
formal and devoid of meaning, begin to take on a life of their own.
For example, it is a standard observation that
combinatorial identities can sometimes be expressed using derivatives of generating
functions.  This is exactly what happens in the case of quantum cohomology,
and it leads to very interesting differential equations which often
have a geometrical meaning.

These notes are a continuation of \cite{Gu}, where
a very brief introduction to the function-theoretic aspects of
quantum cohomology was given.  Needless to say, neither \cite{Gu}
nor the present article is written for experts.  They are meant
to be readable by mathematicians approaching the subject for the first time,
preferably with some background knowledge of differential geometry
and algebraic topology. Much of our exposition is lifted directly from 
\cite{Gi1}-\cite{Gi6} and \cite{Co-Ka}; the only novel aspect lies
in what we have chosen to delete rather than what we have added
(though we have worked out some very detailed examples in the appendices).
Apart from the general goal of explaining how differential
equations arise in quantum cohomology, one of our aims is to arrive
at a point of contact with the \ll mirror phenomenon\rrr.  Even with
such an imprecise goal, our discussion is very incomplete, and in particular we
have to admit that we have not followed up on the loose ends from \cite{Gu}.

Fortunately
there are a number of excellent sources for further information.  The most
comprehensive and the most elementary
is the book \cite{Co-Ka}. The sheer quantity of material
as well as the emphasis on algebraic geometry may be forbidding,
but the careful and helpful presentation makes it invaluable; in addition
there is a lot of new material (in particular new proofs of known results,
and new examples).  Many of the geometric
aspects of the subject have been developed in a series
of fundamental papers by B. Dubrovin, and \cite{Du2} in particular is 
now a classic reference --- although here the tremendous breadth of the
subject matter is a barrier for the beginner.
The papers \cite{Gi1}-\cite{Gi6} of A. Givental provide the
inspiration for the whole subject of \ll homological geometry\rrr, although
the groundbreaking nature of the
arguments makes them hard to follow.  Nevertheless, even for the beginner, Givental's
articles are recommended  because of the wealth of motivation
provided. Other foundational papers are those of 
M. Kontsevich and  Y. Manin
and (for the underlying mathematical physics) those of
E. Witten and C. Vafa.
For a historical perspective and
many more references, \cite{Co-Ka} should be consulted, as well as the recent 
foundational book \cite{Ma}.

We shall generally use the notation of \cite{Gu}; a brief
review follows.   
Let $M$ be a simply connected (and compact, connected) K\"ahler manifold,
of complex dimension $n$, whose nonzero integral cohomology
groups are of even degree and torsion-free.
We choose a basis $b_0,b_1,\dots,b_s$ of $H^\ast(M;\Z)$,
such that $b_1,\dots,b_r$ form a basis of $H^2(M;\Z)$.  The
Poincar\acuteaccent e dual basis of $H_\ast(M;\Z)$ will be denoted
by $B_0,B_1,\dots,B_s$.  The dual basis of $H^\ast(M;\Z)$ with respect to
the intersection form $(\ ,\ )$ will be denoted
by $a_0,a_1,\dots,a_s$.  Thus, we have
$(a_i,b_j)=\lan a_i,B_j\ran = \lan b_i, A_j \ran = \de_{ij}$.
We shall choose $b_0=1$, the identity element of $H^\ast(M;\Z)$,
so that $B_0$ is the fundamental homology class of $M$; sometimes we
write $B_0=M$.  

For homology classes (or representative
cycles of such classes --- we blur the distinction) $X_1,\dots,X_i$, when $i\ge 3$,
the notation $\lan X_1 \vert \dots \vert X_i
\ran\subD$ will denote
the \ll usual\rr genus $0$ Gromov-Witten invariant
obtained using moduli spaces of \ll stable rational curves with $i$ marked
points\rrr. For the definition and properties of 
$\lan X_1 \vert \dots \vert X_i\ran\subD$ we refer to \cite{Fu-Pa}
and chapter 7 of \cite{Co-Ka} (where the standard notation
$\langle I_{0,i,D}\rangle(x_1,\dots,x_i)_{0,D}$ is used).
Here, $D$ is an element of $\pi_2(M)$, so we may write $D=\sum_{i=1}^r s_i A_i$,
and we shall assume as in \cite{Gu} that the homotopy class $D$ contains
holomorphic maps $\C P^1\to M$ only when $s_i\ge 0$ for all $i$.

It is necessary to issue a warning at this point.  In \S 7 of \cite{Gu}, the notation 
$\lan X_1 \vert \dots \vert X_i\ran\subD$ had a different meaning, 
namely the intersection number
$\vert
\HolD^{X_1,p_1} \cap \dots \cap \HolD^{X_i,p_i} 
\vert$, where the points $p_1,\dots,p_i$ are fixed. 
To avoid confusion the latter will be denoted by
$\lan X_1 \vert \dots \vert X_i \ran^{fix}\subD$ in the present article.
For $i=3$, the two definitions agree. For $i\ge 4$, they are
(in the words of \cite{Fu-Pa}) solutions to two different
enumerative problems, and they have somewhat different properties.

A general element of $H^\ast(M;\C)$  will be denoted by
$\hat t=\sum_{i=0}^s t_i b_i$. Since elements of $H^2(M;\C)$ play
a special role, we reserve the symbol $t$ for a general element of $H^2(M;\C)$,
i.e. $t=\sum_{i=1}^r t_i b_i$.

The \ll large\rr quantum product on the vector space $H^\ast(M;\C)$
is defined by
$$
\lan a\cth b,C\ran = \sum_{D,k\ge 0} \frac{1}{k!}
\lan A\vert B \vert C \vert \hat T  \lkr  \hat T \ran\subD
$$
where (as the notation indicates)
the Poincar\acuteaccent e dual $\hat T$ of $\hat t$ appears
$k$ times in the general term of the series. The $D=0$ term is special, because 
$\lan A\vert B \vert C \vert \hat T  \lkr  \hat T \ran_0$ is zero unless
$k=0$.  Using this fact, and the \ll divisor rule\rr
$\lan A \vert B \vert C \vert T \lkr T \ran\subD=
\lan A \vert B \vert C  \ran\subD \lan t, D\ran^k$, we see that
the \ll small\rr quantum product $a\ct b$ is equal to the quantum product which
was used in \cite{Gu}.
We shall not be concerned with the question of convergence
of infinite series like this; we shall assume that the series
converges in a suitable region or simply treat it as a formal series.
Each of $\cth$ and $\ct$ endows $H^\ast(M;\C)$ with the structure of a
commutative algebra (over $\C$) with identity element $b_0=1$.

The (large) quantum product on the vector space $H^\ast(M;\C)$ is
determined by giving all quantum products of the basis elements $b_i$;
these in turn are determined by the following function,
which is called the Gromov-Witten potential:
$$
\Phi(\hat t)= \sum_{D,k\ge 3} \frac{1}{k!}
\lan \hat T  \lkr  \hat T \ran\subD.
$$
This may be regarded as a generating function for the
Gromov-Witten invariants; it is rather unwieldy, of course, and
one of the main themes of the subject is the fact that there
are alternative expressions for it.  
Because of the linearity of Gromov-Witten invariants, we have
$$
\ddi \lan \hat T \vert A \vert B \vert \dots \ran = 
 \lan B_i \vert A \vert B \vert \dots \ran
$$
and hence
$$
\lan b_i\cth b_j, B_k \ran = \frac {\b^3} {\b t_i \b t_j \b t_k } \Phi 
= b_i b_j b_k \Phi
$$
(we identify $b_i$ with the constant vector field $\ddi$ when it
is convenient to do so). This is how derivatives of generating
functions enter into the theory.

In \S 1 we review the definition of the Dubrovin connection, and
its \ll fundamental matrix\rr $H$ of flat sections. This serves as a
generating function for a special family of 
Gromov-Witten invariants, as we explain in \S 2, and it is this which we shall
use rather than $\Phi$.  The fact that it, like $\Phi$, 
has an alternative characterization, as a solution of the
\ll quantum differential equations\rrr, is our main focus, and we shall illustrate it
by means of various examples (in \S 3 and in the
appendices).   In \S 4 we discuss briefly the extent to which these examples can
be generalized.  

In all such cases, the quantum differential equations can be
solved explicitly in terms of generalized {\it hypergeometric functions.}
This leads to the great surprise of the subject: the unruly Gromov-Witten
invariants are governed (in these examples) by a very simple underlying
principle, in the sense that their generating functions are given by
explicit formulae.  In particular, it is easy to compute Gromov-Witten
invariants this way, whereas computations using the original
definition overwhelm all known techniques of algebraic geometry
almost immediately.  This phenomenon, conjectured by physicists,
was the original motivation for homological geometry.

For mathematicians, there are two serious difficulties here (and we do not claim to
shed any light on them in this article).  The first is that it is desirable to
prove {\it a priori} that the generating functions are given by specific
hypergeometric functions for large classes of manifolds.  This is the
\ll Mirror Theorem\rr or \ll Mirror Identity\rrr, and it has in fact been proved in various
situations, in particular in the cases originally discussed by physicists.  
The second is the question of {\it why} such results should be true. This
is related to the \ll Mirror Symmetry Conjecture\rrr, for which
there is as yet no mathematical foundation.

{\eightpoint
Acknowledgements:  The author is grateful to Augustin-Liviu Mare and
Takashi Otofuji for several helpful discussions.  Part of this work
was done in September 2000 while the author was visiting the 
UNAM in Cuernavaca and the CINVESTAV in Mexico City,
and he is grateful to both institutions --- and
in particular to Jose Seade of UNAM, and Luis Astey and Elias Micha
of CINVESTAV --- for their hospitality.}

\head
\S 1 The Dubrovin connection
\endhead

References:  \cite{Du1}, \cite{Du2}, \cite{Ma}, \cite{Co-Ka}, \cite{Gu}

The formula $\om_{\hat t}(x)(y)=x\cth y$ defines a $1$-form on (the complex
manifold) $W=H^\ast(M;\C)$ with values in the Lie algebra
of endomorphisms of (the complex vector space) $W$.  The Dubrovin connection on 
$W$ is defined by $\nabla^{\la} = d + \la \om$. The
properties of the quantum product (see \cite{Gu}) give:

\proclaim{Theorem} The Dubrovin connection has zero curvature (i.e.
is flat) for every value of $\la\in\C$.  In other words,
$d\om=\om\wedge\om=0$.
\endproclaim

\no It follows (from $d\om=0$ and $d\om +\la \om\wedge\om=0$, respectively) 
that there exist functions

$K:W\to \End(W)$ such that $\la\om = dK$, and

$H:W\to Gl(W)$ such that $\la\om = dH H^{-1}$.

\no By elementary properties of first-order differential
equations, each function is determined uniquely when
its value is specified at a single point of $W$.

It is easy to verify that a suitable function $K$ is given explicitly by
$$
\align
\lan K(\hat t)(a),C \ran &=
\la\sum_{D\ge 0,k\ge 1} \frac {1}{k!} 
\lan A\vert C\vert \hat T  \lkr  \hat T \ran\subD\\
&=\la\lan A\vert C\vert \hat T \ran_0 \plus
\la\sum_{D\ne 0,k\ge 1} \frac {1}{k!} 
\lan A\vert C\vert \hat T  \lkr  \hat T \ran\subD.
\endalign
$$  
Alternatively, by making use of the chosen bases, we may write this formula as
$$
K(\hat t)(a) =
\la\sum_{j=0,\dots,s}\lan A\vert B_j\vert \hat T \ran_0 \ a_j \plus
\la\sum_{D\ne 0,k\ge 1,j=0,\dots,s} \frac {1}{k!} 
\lan A\vert B_j\vert \hat T  \lkr \hat T \ran\subD \ a_j.
$$
Replacing $k\ge 1$ by $k\ge 0$ in these formulae adds a constant to
$K$, and hence produces another solution.  (The Gromov-Witten invariant
$\lan X_1 \vert \dots \vert X_i\ran\subD$ can be defined for any
$i\ge 0$, so the formulae make sense.)
The restriction of this modified $K$ to $H^2(M;\C)$ simplifies
(on making use of the divisor rule) to
$$
K( t )(a)  = 
\la\sum_{j=0,\dots,s} \lan A\vert B_j\vert T\ran_0 \ a_j +
\la\sum_{D\ne 0, j=0,\dots,s} \lan A\vert B_j\ran\subD e^{\lan t,D\ran} \ a_j,
$$
which is the formula in \S 7 of \cite{Gu}.

Our main interest will be the function $H$, however, for which an
explicit formula is less obvious.  From now on,
following Givental, we shall make the change of notation
$
\la =  1/h
$
so that the equation $dH H^{-1}=\la\om$ becomes
$hdH = \om H$.  More explicitly, for $\hat t,x\in H^\ast(M;\C)$
this equation is
$
h (dH)_{\hat t} (x) = \om_{\hat t}(x)H(\hat t).
$
The corresponding equations for the column vectors
$
H_i(\hat t) = H(\hat t)(b_i), 0\le i \le s,
$
of $H$ are
$$
h \ddj H_i = b_j\cth H_i,\quad 0\le i,j \le s.
$$
We can write the original equations 
for the matrix-valued function $H:W\to \End(W)$ in the form
$$
h \ddj H = M_j H,\quad 0\le j \le s,
$$
where $M_j$ denotes the matrix of the quantum multiplication
operator $b_j\cth$.  We could write these equations
more informally as $h \ddj H = b_j\cth H$, with
the understanding that $b_j\cth$ operates on the column
vectors of $H$.

Alternatively, we may consider the equations
$$
h \ddj \psi = b_j\cth \psi,\quad 0\le j \le s
$$
for a vector-valued function $\psi:W\to W$.  The solution
space of this system is (under favourable conditions) $s+1$-dimensional, and 
any basis $\psi_{(0)},\dots,\psi_{(s)}$ of
solutions gives rise to a matrix-valued function 
$$
H=
\pmatrix
\vert & & \vert \\
\psi_{(0)} & \dots & \psi_{(s)}\\
\vert & & \vert
\endpmatrix
$$ 
of the required type.

The case $j=0$ is particularly simple, as $b_0\cth x = x$ for any 
cohomology class $x$ (we have chosen $b_0=1$, and this is 
the identity element in the quantum cohomology algebra as well as
in the ordinary cohomology algebra).  Hence the $t_0$-dependence of $H$ or $\psi$ is
just given by a factor of $e^{t_0/h}$.

The following precise connection between $H$ and the quantum product was
established in \cite{Gi-Ki}.  We restrict attention to $t\in H^2(M;\C)$ as this
is the version that will be needed later.

\proclaim{Theorem}  Let $P(X_0,\dots,X_{2r})$ be a polynomial in
$2r+1$ variables, written so that, in each monomial term, $X_{i}$
precedes $X_{j}$ if $1\le i\le r$ and $r+1\le j\le 2r$. If the
functions $(H_u,1)$, $u=0,\dots,s$, satisfy the differential equation
$
P(h,e^{t_1},\dots,e^{t_r},h\ddone,\dots,h\ddare)(H_u,1)=0,
$
then the relation $P(0,e^{t_1},\dots,e^{t_r}, b_1,\dots,b_r)=0$
holds in the quantum cohomology algebra $( H^\ast(M;\C),\ct)$.
\endproclaim

\demo{Proof} This is a simple consequence of the relation
$h \ddj H_i = b_j\ct H_i$, although some care is needed because
of the fact that the quantum product depends on $t$.   The parameter
$h$ will play a role in the argument.
First, note that
$$
\align
h\ddi (H_u,f) &= (h\ddi H_u, f) + (H_u, h\ddi f) \\
&= (b_i\ct H_u, f) + (H_u, h\ddi f) \\
&= (H_u, b_i \ct f) + (H_u, h\ddi f) \\
&= (H_u, b_i \ct f + h\ddi f )
\endalign
$$
for any function $f:H^2(M;\C) \to H^\ast(M; \C)$ and for any
$u=0,\dots,s$.  Repeated application of this formula shows that
$$
\align
P(h,e^{t_1},\dots,e^{t_r},&h\ddone,\dots,h\ddare)(H_u,f)=\\
&(H_u, 
P(h,e^{t_1},\dots,e^{t_r}, b_1 \ct  + h\ddone ,\dots,b_r \ct  + h\ddare)f ).
\endalign
$$
The operators $b_i\ct$ and $h\ddj$ do not necessarily commute here.

Take $f$ to be the constant function $1$ in this formula.  
Then by hypothesis the left hand side  is zero for all $u$, hence
the right hand side is zero for all $u$, and so (as $H_0(t),\dots,H_s(t)$ are a basis) 
$P(h,e^{t_1},\dots,e^{t_r}, b_1 \ct  + h\ddone ,\dots,b_r \ct  + h\ddare)1=0$.
This holds for all values of $h$ ($\ne 0$), so we obtain (as $h\to 0$)
$P(0,e^{t_1},\dots,e^{t_r}, b_1,\dots,b_r)=0$, as required.
\qed
\enddemo

As special cases of the formula in the proof, we have
$$
\align
h\ddi (H_u,1) &= (H_u, b_i\ct 1) = (H_u,b_i)\\
(h\ddj)(h\ddi)(H_u,1) &= (H_u, b_j\ct b_i + h\ddj b_i ) = (H_u,b_j\ct b_i )
\endalign
$$
from which it follows that the converse of the theorem is true when $P$
is of degree at most two.

\head
\S 2 $H$ via descendent Gromov-Witten invariants
\endhead

References:  \cite{Gi1}-\cite{Gi6}, \cite{BCPP}, \cite{Co-Ka}, \cite{Pa}

Remarkably, there is an explicit formula for $H$,
which uses {\it equivariant} Gromov-Witten invariants.  By using a
localization theorem, this formula may then be expressed in terms
of {\it descendent} Gromov-Witten invariants.  We shall discuss the
latter formula in this section. Full details can be found in
chapter 10 of \cite{Co-Ka}. 

The descendent Gromov-Witten invariants
$$
\lan \tau_{d_1} X_1 \vert  \tau_{d_2} X_2 \vert \dots \vert  \tau_{d_i} X_i \ran\subD
$$
(which originate from the
gravitational descendents or gravitational correlators of physics)
are generalizations
of the {\it primary} Gromov-Witten invariants
$\lan X_1 \vert X_2 \vert \dots \vert X_i \ran\subD$, to which they
specialize when $d_1=d_2=\dots=d_i=0$.  Whereas the primary invariants
are defined using cycles in the moduli space which come from
cycles $X_1,X_2,\dots,X_i$ in $M$, the descendent invariants incorporate
additional cycles representing the first chern classes of certain
line bundles $\Cal L_1^{d_1},\Cal L_2^{d_2},\dots \Cal L_i^{d_i}$ on the moduli space.
(The $d_1,d_2,\dots,d_i$ are nonnegative integers.) 
A necessary condition for 
$\lan \tau_{d_1} X_1 \vert  \tau_{d_2} X_2 \vert \dots \vert  \tau_{d_i} X_i \ran\subD\ne 0$
is the numerical condition
$$
\sum_{j=1}^i \vert x_j\vert  \plus  2\sum_{j=1}^i  d_j
= 2(n+i-3) \plus 2\lan c_1TM,D\ran
$$
where $n$ is the complex dimension of $M$.
(This is called the \ll degree axiom\rrr.)

\proclaim{Theorem} A solution $H:W\to \End(W)$ of the system
$h \ddj H = b_j\cth H$, $0\le j\le s$, is given
by 
$$
H(\hat t)(a)=a\plus
\sum_{D\ge 0,k\ge 0,n\ge 0,j=0,\dots,s}
\frac 1 {h^{n+1}}
\frac 1 {k!}
\lan \tau_nA \vert B_j \vert \hat T  \lkr  \hat T \ran\subD \ a_j
$$
where $\hat t,a\in W$.
\endproclaim

\demo{Proof} Let
$$
\llan X_1\vert \dots \vert X_i\rran = \sum_{D\ge 0,k\ge 0} \frac 1{k!} 
\lan X_1\vert \dots \vert
X_i\vert \hat T \lkr  \hat T\ran\subD.
$$
With this (standard) notation the definition of the quantum product becomes
$$
\lan a\cth b, C\ran = \llan A \vert B \vert C \rran,
$$
which is quite analogous to the definition of the ordinary cup product,
although one must remember that the variable $\hat t$ is implicit in $\llan\ ,\ \rran$.
The above formula for $H$ becomes
$$
H(\hat t)(a) = a \plus 
\sum_{n\ge 0,j=0,\dots,s} \frac 1{h^{n+1}} \llan \tau_n A\vert B_j \rran a_j
$$
and so we must prove that this function satisfies the equations 
$h\ddi H(\hat t)(a) = b_i \cth H(\hat t)(a)$ for $i=0,\dots,s$.

We have $\ddi \lan \tau_n A \vert C \vert \hat T 
\lkr \hat T\ran\subD =
k\lan \tau_n A \vert C \vert B_i \vert \hat T  
\lkmor  \hat T \ran\subD$, from which it follows that
$\ddi \llan \tau_n A \vert C \rran = \llan \tau_n A \vert C \vert B_i \rran$.
Hence the left hand side of the equation is
$$
\align
h\ddi H(\hat t)(a) &= \sum_{n\ge 0,j=0,\dots,s} \frac 1 {h^n}
\llan \tau_n A \vert B_j \vert B_i \rran a_j\\
&=\sum_{j=0,\dots,s}
\llan \tau_0 A \vert B_j \vert B_i \rran a_j
\plus
\sum_{n\ge 1,j=0,\dots,s} \frac 1 {h^n}
\llan \tau_n A \vert B_j \vert B_i \rran a_j\\
&=
\sum_{j=0,\dots,s}
\llan A \vert B_j \vert B_i \rran a_j
\plus
\sum_{n\ge 0,j=0,\dots,s} \frac 1 {h^{n+1}}
\llan \tau_{n+1} A \vert B_j \vert B_i \rran a_j.
\endalign
$$
The right hand side is
$$
\align
b_i \cth H(\hat t)(a) &=
b_i\cth a + \sum_{n\ge 0,j=0,\dots,s} \frac 1 {h^{n+1}}
\llan \tau_n A \vert B_j\rran b_i\cth a_j\\
&= \sum_{j=0,\dots,s}
\llan A \vert B_i \vert B_j \rran a_j
\plus
\sum_{n\ge 0,j,u=0,\dots,s} \frac 1 {h^{n+1}}
\llan \tau_n A \vert B_j\rran \llan B_i \vert A_j \vert B_u\rran a_u.
\endalign
$$

The first parts of the left and right hand sides agree; to see that
the second parts agree we use the \ll topological recursion relation\rr
$$
\llan \tau_{n+1} A \vert B \vert C \rran =
\sum_{j=0,\dots,s} \llan \tau_n A \vert B_j \rran \llan B \vert C \vert A_j \rran.
$$
This completes the proof.
\qed
\enddemo

By using the divisor property for descendent Gromov-Witten invariants, the
restriction
$
H\vert_{H^2(M;\C)}
$
may be simplified further:

\proclaim{Corollary} For $t\in H^2(M;\C)$ and $a\in W$, a solution of
the system $h \ddj H = b_j\ct H$, $1\le j\le r$, is given by
$$
H(t)(a)= a e^{t/h} \plus
\sum_{D\ne 0,n\ge 0,j=0,\dots,s}
\frac 1 {h^{n+1}}
\lan \tau_n e^{T/h}A \vert B_j \ran\subD \ e^{\lan t,D\ran} \ a_j.
$$
\endproclaim

\no (The notation $T^k A$ means the Poincar\acuteaccent e dual homology
class to $t^ka$.  Extending this convention in an obvious way, it is
convenient to write $\sum_{k\ge 0} T^kA/(h^k k!)=e^{T/h}A$.)

\demo{Proof} Restricting to $t\in H^2(M;\C)$ allows us to make use of
the following divisor rule for descendent Gromov-Witten invariants:

\proclaim{Lemma} Assume that $x\in H^2(M;\C)$. Then
$\lan \tau_{d_1} X_1 \vert \dots \vert \tau_{d_i} X_i \vert X \ran\subD =$
\newline
$
\lan x,D \ran \lan \tau_{d_1} X_1 \vert \dots \vert \tau_{d_i} X_i \ran\subD
+
\sum_{j=1}^{i} 
\lan \tau_{d_1} X_1 \vert \dots \vert \tau_{d_j-1} XX_j 
\vert \dots \vert \tau_{d_i} X_i \ran\subD.
$
\qed
\endproclaim

Applying this, and writing
$\al_{n,k,l}=\lan \tau_nA T^l \vert C \vert T  \lkr  T \ran\subD$,
we obtain
$$
\align
\lan \tau_nA \vert C \vert T  \lkr   T \ran\subD &= \al_{n,k,0}\\
&= \lan t,D \ran \al_{n,k-1,0} + \al_{n-1,k-1,1}\quad\text{(the divisor rule)}\\
&= \quad\quad\cdots\\
&= \sum_{\mu=0}^{k} \tfrac {k!} {\mu! (k-\mu)!}  \lan t,D \ran^{\mu}
\al_{n-(k-\mu),0,k-\mu}\\
&=\sum_{\mu=0}^{k} \tfrac {k!} {\mu! (k-\mu)!}  \lan t,D \ran^{\mu}
\lan \tau_{n-(k-\mu)} AT^{k-\mu} \vert C \ran\subD.
\endalign
$$
In the formula of the theorem for $H(t)$, we consider first the
term with $D=0$.  This becomes
$$
\align
\sum_{k\ge 0,n\ge 0,j=0,\dots,s}
&\frac 1 {h^{n+1}}
\frac 1 {k!}
\lan \tau_nA \vert B_j \vert T  \lkr  T \ran_0 \ a_j\\
&=
\sum_{k\ge 0,n\ge 0,j=0,\dots,s}
\frac 1 {h^{n+1}}
\frac 1 {k!}
\lan \tau_{n-k+1}  AT^{k-1} \vert B_j \vert T \ran_0 \ a_j.
\endalign
$$
Now, it can be shown that
$\lan \tau_{n-k+1} AT^{k-1} \vert B_j \vert T \ran_0$
can be nonzero only for $k=n+1$, in which case it is 
$\lan  AT^{n} \vert B_j \vert T \ran_0 = \lan ab_j t^{n+1},M\ran$.
Hence the sum is
$$
\sum_{j=0,\dots,s} \left( \frac 1h \lan ab_jt,M\ran + 
\frac 1{h^2} \frac 1{2!} \lan ab_jt^2,M\ran +
\frac 1{h^3} \frac 1{3!} \lan ab_jt^3,M\ran +    \dots\right) a_j.
$$
If we add the first term \ll$a$\rr in the formula of the theorem,
in the form $a=\sum_{j=0,\dots,s} \lan a,B_j\ran a_j
= \sum_{j=0,\dots,s} \lan ab_j,M\ran a_j$,
we obtain $a e^{t/h}$.

Next we consider the terms with $D\ne 0$. These are
$$
\align
\sum_{D\ne 0,k\ge 0,n\ge 0,j=0,\dots,s}
&\frac 1 {h^{n+1}}
\frac 1 {k!}
\lan \tau_nA \vert B_j \vert \hat T  \lkr  \hat T \ran\subD \ a_j\\
&=\sum_{D\ne 0,k\ge 0,n\ge 0,j=0,\dots,s,\mu+\nu=k}
\frac 1 {h^{n+1}}
\frac 1{\mu!\nu!}
\lan t,D\ran^{\mu} \lan \tau_{n-\nu} AT^{\nu} \vert B_j \ran\subD \ a_j\\
&=
\sum_{D\ne 0,n\ge 0,j=0,\dots,s}
\quad
\sum_{\mu\ge 0,\nu\ge 0}
\frac 1 {h^{n+1}}
\frac 1{\mu!\nu!}
\lan t,D\ran^{\mu} \lan \tau_{n-\nu} AT^{\nu} \vert B_j \ran\subD \ a_j\\
&=
\sum_{D\ne 0,n\ge 0,j=0,\dots,s}
\quad
\sum_{\nu\ge 0}
e^{\lan t,D \ran}
\frac 1 {h^{n+1}}
\frac 1{\nu!}
\lan \tau_{n-\nu} AT^{\nu} \vert B_j \ran\subD \ a_j\\
&=
\sum_{D\ne 0,j=0,\dots,s}
e^{\lan t,D \ran}
\sum_{n\ge\nu\ge 0}
\frac 1 {h^{n+1}}
\frac 1{\nu!}
\lan \tau_{n-\nu} AT^{\nu} \vert B_j \ran\subD \ a_j\\
&=
\sum_{D\ne 0,j=0,\dots,s}
e^{\lan t,D \ran}
\sum_{\nu\ge 0,l\ge 0}
\frac 1 {h^{l+1}}\frac 1 {h^{\nu}}
\frac 1{\nu!}
\lan \tau_{l} AT^{\nu} \vert B_j \ran\subD \ a_j
\quad\text{($l=n-\nu$)}\\
&=
\sum_{D\ne 0,l\ge 0,j=0,\dots,s}e^{\lan t,D \ran}
\frac 1{h^{l+1}} 
\lan \tau_{l}Ae^{T/h} \vert B_j \ran\subD \ a_j
\endalign
$$
This completes the proof.
\qed
\enddemo

\head
\S 3 The quantum differential equations
\endhead

References:  \cite{Gi1}-\cite{Gi6}, \cite{BCPP}, \cite{Co-Ka}, \cite{Pa}

We are now ready to describe the crucial phenomenon.  So far we have
seen that there are various natural generating functions for
(various families of) Gromov-Witten invariants. It turns out
that they satisfy meaningful differential
equations.  To establish this phenomenon rigorously for 
general classes of manifolds is a difficult problem 
which is an area of current research.   We shall merely
give statements and some  illustrative examples.   The reader should bear
in mind that, although our examples are important, they are nevertheless
very special.

The generating functions that we shall focus on are the rows of
the matrix function $H(t)$, so
we shall use matrix notation more systematically in this section. All
matrices are taken with respect to the basis $b_0,\dots,b_s$. The product
of two matrices $M_1,M_2$ will be denoted by $M_1\cdot M_2$.  When we use
this notation we shall regard
$b_0,\dots,b_s$ as column vectors and $a_0,\dots,a_s$ as row vectors.
Hence we have $a_i\cdot b_j = (a_i,b_j) = \de_{ij}$, and
$$
H=
\pmatrix
\vert & & \vert \\
H_0 & \dots & H_s\\
\vert & & \vert
\endpmatrix
=
\pmatrix
- & J_0 & - \\
  &  \vdots & \\
- & J_s & - 
\endpmatrix
$$
where the columns and rows of $H$ are, respectively,
$$\align
H_i(t)&=H(t)(b_i)=H(t)\cdot b_i= \sum_{j=0,\dots,s} H_{ji}(t) b_j\\ 
J_i(t)&=a_i\cdot H(t)= \sum_{j=0,\dots,s} H_{ij}(t) a_j.
\endalign
$$
The $(i,j)$-th entry of $H$ is 
$H_{ij}=J_i\cdot b_j = a_i\cdot H_j = a_i\cdot H\cdot b_j$.

Without loss of generality we assume that $a_s=b_0=1$, so the
last row of $H$ is simply 
$$
J_s(t)=\sum_{j=0,\dots,s} H_{sj}(t)a_j 
=\sum_{j=0,\dots,s} a_s\cdot H(t)\cdot b_j \ a_j
=\sum_{j=0,\dots,s} (a_s, H_j(t) )a_j
=\sum_{j=0,\dots,s} (1, H_j(t)) a_j.
$$
Following Givental we shall denote
this function by $J$ from now on:

\proclaim{Definition} 
$J(t)=J_s(t)=\sum_{j=0,\dots,s} (H_j(t),1) a_j.$
\endproclaim

\no One also has the more canonical expression $J(t)=H(t)^\dagger (1)$,
where $H(t)^\dagger$ is the adjoint of $H(t)$ with respect to the
bilinear form $(\ ,\ )$.  The theorem of \S 1
relating
$H$ with the quantum product can be re-stated as follows: 

\proclaim{Theorem}  Let $P(X_0,\dots,X_{2r})$ be a polynomial in
$2r+1$ variables, written so that, in each monomial term, $X_{i}$
precedes $X_{j}$ if $1\le i\le r$ and $r+1\le j\le 2r$. If the
function $J$ satisfies the differential equation
$
P(h,e^{t_1},\dots,e^{t_r},h\ddone,\dots,h\ddare)J=0,
$
then the relation $P(0,e^{t_1},\dots,e^{t_r}, b_1,\dots,b_r)=0$
holds in the quantum cohomology algebra $( H^\ast(M;\C),\ct)$.
\endproclaim

For this reason we shall investigate the function $J$ in more detail.
It, and more generally the $i$-th row $J_i$, can be written explicitly in
terms of Gromov-Witten invariants as follows:

\proclaim{Proposition}
$$
J_i(t) = e^{t/h} \left( a_i \plus 
\sum_{D\ne 0,n\ge 0,j=0,\dots,s} \frac 1 {h^{n+1}} \lan \tau_n B_j \vert A_i\ran\subD
\ e^{\lan t,D\ran} \ a_j \right).
$$
\endproclaim

\demo{Proof} Substituting the earlier formula for $H_j(t)$ 
in the definition of $J_i(t)$, we obtain:
$$
\align
J_i(t) &= \sum_{j=0,\dots,s} (H_j(t),a_i) a_j \\
&=
\sum_{j=0,\dots,s} \left(
b_j e^{t/h} \plus
\sum_{D\ne 0,n\ge 0,u=0,\dots,s}
\frac 1 {h^{n+1}}
\lan \tau_n e^{T/h}B_j \vert B_u \ran\subD \ e^{\lan t,D\ran} \ a_u,a_i
\right) 
a_j\\
&=\sum_{j=0,\dots,s} \left(
b_j e^{t/h} \plus
\sum_{D\ne 0,n\ge 0,u=0,\dots,s}
\frac 1 {h^{n+1}}
\lan \tau_n e^{T/h}B_j \vert A_u \ran\subD \ e^{\lan t,D\ran} \ b_u,a_i
\right) 
a_j\\
&=
\sum_{j=0,\dots,s} (a_i e^{t/h},b_j)a_j \plus
\sum_{D\ne 0,n\ge 0,j=0,\dots,s}
\frac 1 {h^{n+1}}
\lan \tau_n e^{T/h}B_j \vert A_i \ran\subD \ e^{\lan t,D\ran} a_j\\
&=a_i e^{t/h} \plus 
\sum_{D\ne 0,n\ge 0,j=0,\dots,s}
\frac 1 {h^{n+1}}
\lan \tau_n e^{T/h}B_j \vert A_i \ran\subD \ e^{\lan t,D\ran} a_j
\endalign
$$
Replacing $b_0,\dots,b_s$ by the basis
$e^{-t/h}b_0,\dots,e^{-t/h}b_s$, we get
$$
\align
J_i(t) &= a_ie^{t/h} \plus 
\sum_{D\ne 0,n\ge 0,j=0,\dots,s} \frac 1 {h^{n+1}} \lan \tau_n B_j \vert A_i\ran\subD
\ e^{\lan t,D\ran} \ e^{t/h} a_j\\
&=e^{t/h} \left( a_i \plus 
\sum_{D\ne 0,n\ge 0,j=0,\dots,s} \frac 1 {h^{n+1}} \lan \tau_n B_j \vert A_i\ran\subD
\ e^{\lan t,D\ran} \ a_j \right)
\endalign
$$
as required.
\qed
\enddemo

\proclaim{Corollary}
$$
J(t) = e^{t/h} \left( 1 \plus 
\sum_{D\ne 0,n\ge 0,j=0,\dots,s} \frac 1 {h^{n+1}} \lan \tau_n B_j \vert M\ran\subD
\ e^{\lan t,D\ran} \ a_j \right).
$$
\endproclaim

In the following concrete examples, we shall use our knowledge of the
quantum product to compute $J$ and $H$ explicitly, in order to illustrate
the general theory.   In  applications, however, one would hope to
do the opposite, i.e. use (at least partial) knowledge of $J$ or $H$
to investigate the quantum product.  This strategy has been
used successfully in \cite{Gi3} and \cite{Ki1}.

\proclaim{Example: $M=\C P^m$}
\endproclaim

The quantum cohomology algebra of $\C P^m$ is well known and
was explained in detail in \S 3 of \cite{Gu}.  We shall
use the same notation here.  Thus we have
$H^\ast(\C P^m;\C) = \oplus_{i=0}^m \C x_i$, and
$x_i=x^i$ where $x=x_1$ is a generator of $H^2(\C P^m;\C)$,
and we write $\hat t = \sum_{i=0}^m t_i x^i$ and $t=t_1 x$.  

With these conventions, the differential equation for
$\psi:H^2(\C P^m;\C) \to H^\ast(\C P^m;\C)$ is
$$
h \ddtone \psi = x \ct \psi.
$$
Let us first note that this can be written in terms
of $\psi= \psi_0  + \psi_1 x + \dots + \psi_m x^m$ as
$$
\align
h(\psi_0^\prime  + \dots + \psi_m^\prime x^m) &=
x\ct(\psi_0  + \dots + \psi_m x^m)\\
&=\psi_0 x + \dots + \psi_{m-1} x^m + e^{t_1} \psi_{m}
\endalign
$$
(where prime denotes derivative with respect to $t_1$), i.e.
$$
h \psi_m^\prime = \psi_{m-1},\quad
h\psi_{m-1}^\prime = \psi_{m-2},\quad \dots, \quad
h\psi_1^\prime = \psi_0,\quad
h\psi_0^\prime = e^{t_1} \psi_m.
$$
Hence $\psi$ may be expressed as
$
\psi = h^m f^{(m)}  + \dots + h f^\prime x^{m-1} + f x^m
$
where $f$ is a solution of the \ll quantum differential equation\rr
$( h\ddtone )^{m+1} f = e^{t_1} f$.

Alternatively, if we write
$$
H=
\pmatrix
- & J_0 & - \\
  &  \vdots & \\
- & J_m & - 
\endpmatrix
$$
then the equation $h \ddtone H = x\ct H$ is 
$$
h
\pmatrix
- & J_0^\prime & - \\
  &  \vdots & \\
- & J_m^\prime & - 
\endpmatrix
=
\pmatrix
0 &  & & e^{t_1} \\
1 & \ddots & & \\
 & \ddots & 0 & \\
 & & 1 & 0
\endpmatrix
\pmatrix
- & J_0 & - \\
  &  \vdots & \\
- & J_m & - 
\endpmatrix
$$
where the first matrix on the right hand side is
the matrix of $\ x\ct\ $ with respect to the basis
$x_0.\dots,x_m$.  This is equivalent to the system
$$
h J_m^\prime = J_{m-1},\quad
hJ_{m-1}^\prime = J_{m-2},\quad \dots, \quad
hJ_1^\prime = J_0,\quad
hJ_0^\prime = e^{t_1} J_m,
$$
which in turn is equivalent to the quantum differential equation
$(h\ddtone)^{m+1}J = e^{t_1}J$ for $J=J_m$, with
$$
H=
\pmatrix
- & h^mJ^{(m)} & - \\
  &  \vdots & \\
- & hJ^\prime & -\\
- & J & - 
\endpmatrix.
$$

One can verify that the $H^\ast(\C P^m;\C)$-valued map
$$
J(t) = f_{(0)}(t) + f_{(1)}(t) + \dots + f_{(m)}(t) x^m
=\sum_{d\ge 0} \frac
{ e^{(x/h + d)t_1}  }
{ [ (x+h)(x+2h)\dots(x+dh) ]^{m+1} }
$$
is a solution.  Indeed, writing $q_1=e^{t_1}$, so that
$\ddtone = q_1 \ddqone$, we have
$$
J(t)
=\sum_{d\ge 0} \frac
{ q_1^{x/h + d}  }
{ [ (x+h)(x+2h)\dots(x+dh) ]^{m+1} }, 
$$
and since
$hq_1\ddqone  q_1^{x/h + d} = (x+dh)  q_1^{x/h + d}$,
it follows immediately that 
$$
(hq_1\ddqone)^{m+1} J = q_1 J.
$$
Note, in particular, that
$$
f_{(0)}(t)=\sum_{d\ge 0} \frac{q_1^d}{(d!)^{m+1}h^{d(m+1)}},
$$
which is a hypergeometric series.   The function $J(t)$ may be regarded
as a {\it cohomology-valued} hypergeometric series.

From the beginning of this section we have another expression for $J$, namely
$$
\align
J(t)
&=
e^{t/h} \left( 1 + \sum_{d\ge 1, n\ge 0, j=0,\dots,m} 
\frac 1{h^{n+1}} \lan \tau_n X^j \vert \C P^m \ran_d e^{t_1d} x^{m-j}\right).
\endalign
$$
We shall verify the hypergeometric formula
in the case $m=1$, by making
use of the following known descendent Gromov-Witten invariants (see \cite{Co-Ka}):
$$
\lan \tau_{2d-1} \text{\,point\,} \vert \C P^1 \ran_d = \frac 1 {(d!)^2},\quad
\lan \tau_{2d} \C P^1 \vert \C P^1 \ran_d = 
\frac {-2}{(d!)^2} (1 + \frac12 + \dots + \frac1d ).
$$
All other invariants of the form
$\lan \tau_{n} \text{\,point\,} \vert \C P^1 \ran_d$,
$\lan \tau_{n} \C P^1 \vert \C P^1 \ran_d$
are zero, because of the numerical condition (the degree axiom) mentioned earlier.

To emphasize the degrees of the cohomology classes,
we shall write $x_0$ for $1$, and $x_1$ for
$x$.  Noting that
$B_0=\C P^1$, $B_1=\text{\,point\,}$, and $a_0=x_1$, $a_1=x_0$,
we have
$$
\align
J(t)&=
e^{t_1x_1/h} \left( x_0 + \sum_{d\ge 1, n\ge 0} 
\frac 1{h^{n+1}} 
\{
\lan \tau_n \C P^1 \vert \C P^1 \ran_d  x_1
\plus
\lan \tau_n \text{\,point\,} \vert \C P^1 \ran_d  x_0
\}e^{t_1d}
\right)
\\
&=
e^{t_1x_1/h} \left(
x_0 \plus
\sum_{d\ge 1}
e^{t_1d} 
\left\{
\frac 1{h^{2d}}  \frac 1{(d!)^2} x_0 \plus
\frac 1{h^{2d+1}} 
\frac {-2}{(d!)^2} (1 + \frac12 + \dots + \frac1d) x_1
\right\}
\right)
\\
&=
e^{t_1x_1/h} \left(
x_0 \plus
\sum_{d\ge 1} 
\frac {e^{t_1d}} {h^{2d}} 
\frac 1 {(d!)^2}
\left(
x_0-\frac {2x_1}{h} (1 + \frac12 + \dots + \frac1d)
\right)
\right)
\endalign
$$
Since $x_1^2=0$, we have
$$
\align
x_0-\frac {2x_1}{h} (1 + \frac12 + \dots + \frac1d)
&=
[x_0+\frac {2x_1}{h} (1 + \frac12 + \dots + \frac1d)]^{-1}\\
&= 
[x_0+\frac {x_1}{h} (1 + \frac12 + \dots + \frac1d)]^{-2}
\endalign
$$
Hence
$$
\align
J(t)&=
e^{t_1x_1/h} \left(
x_0 \plus
\sum_{d\ge 1} 
\frac{e^{t_1d}}
{[h^d d! (x_0+\frac {x_1}{h} (1 + \frac12 + \dots + \frac1d))]^2}\right)\\
&=e^{t_1x_1/h} 
\sum_{d\ge 0} 
\frac{e^{t_1d}}
{[d! h^dx_0 + d! h^{d-1}(1+\frac12 + \dots + \frac1d)x_1]^2}.
\endalign
$$
This is the same as
$$
e^{t_1x_1/h}
\sum_{d\ge 0} 
\frac{e^{t_1d}}
{[(x_1+hx_0)(x_1+2hx_0)\dots(x_1+dhx_0)]^2}
$$
(the linear term in $x_0,x_1$ is the only nonzero part of the
expansion of the denominator).

Conversely, starting with this formula, one could deduce the formulae
for the descendent Gromov-Witten invariants given earlier.

In Appendices 1 and 2 we shall work out the quantum differential equations and
their solutions for the flag manifold $F_3$ and the Hirzebruch surface
$\Si_1$.  In each case the results are analogous to those for $\C P^n$.  
The main feature of these examples is that the quantum cohomology algebra
is equivalent to a function-theoretic object, as advertised in the
introduction.  We shall just summarize the results of the calculations
here, postponing speculation concerning their meaning or
possible generalization.

First, the quantum product is equivalent (by definition) to the matrix-valued $1$-form $\om$,
and via the Dubrovin connection this is equivalent to the fundamental
matrix of flat sections, $H$.  The linear system of differential equations
defining $H$, namely $h \bi H = M_i H$, may be rewritten as a system of
higher order differential equations for the function $J$, where
$$
H=
\pmatrix
- & J_0 & - \\
  &  \vdots & \\
- & J_s & - 
\endpmatrix
$$
and $J=J_s$ is the last row of $H$.  These higher order equations for
$J$ are the quantum differential equations.  (In the literature, the
quantum differential operators are defined as the  differential
operators which annihilate $J$ --- in our examples, at least, the two concepts
agree.) 

Next, by the process described in \S 1, the quantum differential equations
produce relations in the cohomology algebra.   In our examples, all relations
are obtained in this way.  Thus we may say that the function $J$, or the
quantum differential equations, are equivalent to the relations of the
quantum cohomology algebra.  

This gives a function-theoretic interpretation
of (the isomorphism type) of the quantum cohomology algebra.  However, the
quantum product itself (or equivalently $\om$, or $H$), requires additional
information beyond $J$.  In practice, this amounts to evaluating a finite number of
quantum products $b_i\ct b_j$ of the cohomology basis elements.  (In the
literature this additional information is sometimes called the Pieri formula,
in analogy with classical Schubert calculus.)  In the present framework it
amounts to giving 
\footnote{This is analogous to the elementary procedure of constructing a
first order system of o.d.e. equivalent to the single o.d.e.
$y^{\prime\prime} + ay^{\prime} + by =0$. One way to do this, 
but by no means the only way, is to write $y_1=y, y_2=y^{\prime}$.}
$J_0,\dots,J_{s-1}$ explicitly as polynomials 
$P_0,\dots,P_{s-1}$ in $h$ and
$\bone,\dots,\br$ applied to $J$; 
such expressions follow from the original
system $h \bi H = M_i H$.  

Let us assume that the ordinary cohomology algebra is generated
by two-dimensional classes (as in our examples).  
It follows from the definition of $H$ that the above polynomials 
are the ones which arise when the 
cohomology basis vectors, given originallly as cup products
of two-dimensional classes, are expressed as {\it quantum} products of two-dimensional
classes.  Thus, these are the finite number of quantum products which
have to be evaluated.  Function-theoretically, this corresponds to
a factorization
$$
H=QH_0
$$
where $H_0$ depends only on the ordinary cohomology algebra and the relations
of the quantum cohomology algebra.  The matrix $Q$ represents the 
\ll Pieri formula\rrr.

Finally we shall comment on the relation between the generating functions
of this section and the functions $v(t,q)$, $V(t,q)$ of \S 6 of \cite{Gu}.
Recall that these were defined as follows:
$$
v(t,q)= \sum_{l\ge 0} \frac 1{l!}t\circ \dots\circ t
$$
where $t\circ \dots\circ t$ is
the quantum product of $l$ copies of $t$, and
$$
V(t,q)= \lan v(t,q), M \ran .
$$
We write $v(t,q)$ rather than $v(t)$ to emphasize that $q$ is
considered as a formal variable, independent of $t$.   Note that we
are using the product $\circ$ rather than $\ct$.  The function $V$
has the property that the differential operators annihilating it
correspond to the relations in the quantum cohomology ring
 --- this holds whenever the quantum cohomology
algebra exists and the ordinary cohomology is generated by two-dimensional
classes (see \S 6 of \cite{Gu}).  In fact the same statement holds for
$v$, which suggests a close relation between $v$ and $J$.

Such a relation may be explained by developing the analogous theory
of the connection form $\om$ and the matrix of flat sections $H$
for the case of the quantum product $\circ$.  The connection form
$\om=\sum_{i=1,\dots,r} M_i dt_i$ remains exactly the same, but the
the equation $hd\tilde H \tilde H^{-1}=\om$ gives rise to a different
matrix $\tilde  H$, because we are regarding $q_1,\dots,q_r$ as
constants.  As before we write
$$
\tilde H=
\pmatrix
- & \tilde J_0 & - \\
  &  \vdots & \\
- & \tilde J_s & - 
\endpmatrix
$$
and this is equivalent to a system of higher order differential
equations for $\tilde J_s$ --- the same system as before, except that
$q_1,\dots,q_r$ are regarded as constants.  

The relation betwen $\tilde J$ ($=\tilde J_s$) and the quantum product
is also the same as before, but the statement (and proof - see \S 1) of this
fact can be simplified because the operators $b_i\circ$ and $h\ddj$
commute.  That is, we have
$$
h\ddi (\tilde H_u,f) = (\tilde H_u, b_i \circ f + h\ddi f )
$$
for any function $f:H^2(M;\C) \to H^\ast(M; \C)\otimes\C[q,q^{-1}]$ and for any
$u=0,\dots,s$.  Hence
$$
\align
P(h,q_1,\dots,q_r,&h\ddone,\dots,h\ddare)(\tilde H_u,f)=\\
&(\tilde H_u, 
P(h,q_1,\dots,q_r, b_1 \circ  + h\ddone ,\dots,b_r \circ  + h\ddare)f ).
\endalign
$$
Taking $f$ to be the constant function $1$ gives:

\proclaim{Theorem}  Let $P(X_0,\dots,X_{2r})$ be a polynomial in
$2r+1$ variables. If the function $\tilde J$ satisfies the differential equation
$
P(h,q_1,\dots,q_r,h\ddone,\dots,h\ddare)\tilde J=0,
$
then the relation $P(0,q_1,\dots,q_r, b_1,\dots,b_r)=0$
holds in $( H^\ast(M;\C)\otimes\C[q,q^{-1}],\circ)$.
\endproclaim

\no The parameter $h$ was needed \ll to separate the non-commuting operators 
$b_i\ct$ and $h\ddj$\rrr, but it plays no essential role here.  The statement
of the theorem would be true if the polynomial 
$P(h,q_1,\dots,q_r,h\ddone,\dots,h\ddare)$
were replaced by a polynomial $P(q_1,\dots,q_r,\ddone,\dots,\ddare)$.

Because the quantum differential equation has
constant coefficients, it is easy to solve. 
In fact, it is easy to solve the original system 
$hd \tilde H \tilde H^{-1}= \om\ (=\sum_{i=1,\dots,r} M_i dt_i)$, 
since the matrices $M_i$ commute
and are independent of $t$.  A solution is
$$
\tilde H = e^{\sum t_i M_i/h}
=e^{\sum (t_i b_i\circ)/h}
=e^{(t\circ)/h}.
$$
Hence 
$$
\align
\tilde J(t)&=\sum_{j=0,\dots,s} (\tilde H_j(t),1) a_j\\
&= \sum_{j=0,\dots,s} (\tilde H(t)b_j,1) a_j\\
&= \sum_{j=0,\dots,s} (e^{(t\circ)/h} b_j,1) a_j\\
&= \sum_{j=0,\dots,s} (e^{(t\circ)/h}1, b_j) a_j\\
&=e^{(t\circ)/h}1\\
&= \sum_{l\ge 0}  \frac 1{h^l}\frac 1{l!}t\circ \dots\circ t.
\endalign
$$
Putting $h=1$ gives the function $v$.

Since $v$ is a generating function for the \ll fixed base point\rr
Gromov Witten invariants $\lan X_1 \vert \dots \vert X_i \ran^{fix}\subD$
(see \S 1), and $J$ is a generating function for certain
descendent Gromov-Witten invariants, the above discussion implies
relations between these objects.   In fact it is well known that
the (genus zero) descendent Gromov-Witten invariants contain no more information
than the primary Gromov-Witten invariants (see \cite{Co-Ka}), so
such relations are not surprising.

Both versions lead to a function-theoretic description of ordinary cohomology
by taking appropriate limits.  For $H$ we must take the asymptotic version,
and for $\tilde H$ we simply set $q_1,\dots,q_r=0$.

In view of the fact that $\tilde H$ appears to 
give a much more elementary function-theoretic
version of quantum cohomology, the reader may wonder what advantage there is
to working with the original $H$.   The answer is that it is $H$
which leads to mirror symmetry --- the hypergeometric nature of $H$ is not
accidental and turns out to have deeper geometrical origins.  We will discuss this
in the next section.

\head
\S 4 Computations for other $M$, and the mirror phenomenon
\endhead

References:  \cite{Gi3}-\cite{Gi5}, \cite{EHX},
\cite{Pa}, \cite{Ki2}, \cite{Co-Ka}

The stringent conditions we have imposed on the manifold $M$  make it difficult
to go beyond the above examples (generalized flag manifolds and Fano
toric manifolds).  However, there is an indirect approach, in which one 
considers 

\no (1) a submanifold $M$ of a generalized flag manifold or Fano toric
manifold $X$, and 

\no (2) a version of quantum cohomology using only those
cohomology classes of $M$ which are restrictions of cohomology classes of $X$.

\no This avoids technical problems in working with the quantum cohomology of
$M$ itself, at the cost of losing some information.

The first and most successful example of this strategy was obtained by taking
$X=\C P^n$ and $M$ a hypersurface in $\C P^n$, and in particular
a quintic hypersurface in $\C P^4$ (which happens to be an example
of a Calabi-Yau variety, as required in certain physical theories).  A detailed treatment
of this story may be found in \cite{Co-Ka}.  

When $M$ is a hypersurface in $\C P^n$, the structure of the classical 
cohomology algebra $H^\ast(M;\C)$ is elucidated by the Lefschetz theorems
(see \cite{Gr-Ha}, chapter 1, section 2).  In particular, $H^\ast(M;\C)$
decomposes additively as the sum of 

\no(i) a subalgebra  of $H^\ast(M;\C)$ obtained by restriction of cohomology classes
from $H^\ast(\C P^n;\C)$, and 

\no(ii) a subgroup of
$H^{n-1}(M;\C)$, generated by  so-called primitive elements. 

\no Roughly speaking, the strategy amounts to working in $\C P^n$ but taking 
quantum products only of elements of (i).

\proclaim{Example: $M=Gr_2(\C^4)$}
\endproclaim

As a very simple example, we consider the Grassmannian $M=Gr_2(\C^4)$,
which may be represented as the hypersurface $z_0z_1+z_2z_3-z_4z_5=0$
in $X=\C P^5$ (the Pl\"ucker embedding).   Since the quantum cohomology of
the homogeneous space $Gr_2(\C^4)$ is already known, 
there is no need to make use of this embedding, but we shall do so
as an illustration of what the indirect approach actually computes.

First let us establish some notation for the quantum cohomology of
$Gr_2(\C^4)$.  
Recall (e.g. from section 4 of \cite{Gu})
that the classical cohomology algebra can be written  as
$$
H^\ast(Gr_2(\C^4);\Z) \cong
\Z[c_1,c_2]/\lan c_1^3 - 2c_1c_2,  c_2^2 - c_1^2c_2 \ran.
$$
We may choose additive generators as follows:
$$
\matrix
\quad H^0\quad &\quad H^2 \quad& \quad H^4\quad & 
\quad H^6 \quad& \quad H^8\quad\\
& & & & \\
1 & c_1 & c_2 & c_1 c_2 & c_2^2=c_1^2c_2 \\
 & & c_1^2 &  & \\
\endmatrix
$$
The remaining cup products are determined by the single relation
$c_1^3=2c_1c_2$.

To give a Schubert description of $H^\ast(Gr_2(\C^4);\Z)$, we
introduce
$$
a=c_1, \quad b=c_2, \quad c=c_1^2-c_2 \ (=s_2), \quad d=c_1c_2,\quad z=c_2^2
$$
and choose the following (slightly modified) additive generators:
$$
\matrix
\quad H^0\quad & \quad H^2\quad &\quad H^4\quad & 
\quad H^6\quad &\quad H^8\quad\\
& & & & \\
1 & a & b & d=ab=ac & z=b^2=c^2=ad \\
 & & c &  & \\
\endmatrix
$$
These are subject to the additional relations
$a^2=b+c$, $bc=0$.   It can be verified that the Poincar\acuteaccent e
dual homology classes to these six generators are represented by Schubert
varieties; in particular the Poincar\acuteaccent e dual of the class $a$
is represented by the set of $2$-planes in $\C^4$ whose intersection with the fixed
$2$-plane $\C^2$ has dimension at least one.  This is a Schubert variety in 
$Gr_2(\C^4)$ of (complex) codimension one. In terms of the Pl\"ucker embedding
it is the result of intersecting $Gr_2(\C^4)$ with a hyperplane in $\C P^5$.

From the above relations, it follows that the cohomology of $H^\ast(Gr_2(\C^4);\C)$
decomposes (in the expected manner) as the sum of

\no(i) the five-dimensional subalgebra generated by $a$, with additive generators
$1$, $a$, $a^2=b+c$, $a^3=2d$, $a^4=2z$, and

\no(ii) the one-dimensional subgroup generated by the primitive class 
$b-c\in H^4(Gr_2(\C^4);\C)$.

By similar calculations to those in Appendices 1 and 2
of \cite{Gu}, one obtains the following basic
products. Statements in square brackets are consequences of earlier statements.

\no(1) $a\circ a=a^2=b+c$
\newline(2) $a\circ b=ab=d$
\newline (3) $a\circ c=ac=d$
\newline[4] $a\circ a\circ a=a^3 =2d$
\newline(5) $b\circ c=e^t $
\newline(6) $b\circ b=b^2=z$
\newline(7) $c\circ c=c^2=z$
\newline[8] $a\circ a\circ b=a^2b+e^t =z+e^t $
\newline[9] $a\circ a\circ c=a^2c+e^t =z+e^t $
\newline[10] $a\circ d=ad+e^t =z+e^t $
\newline[11] $a\circ a\circ a\circ a = a^4 + 2e^t =2z+2e^t $

Of these, (1)-(3) and (5)-(7) determine all possible
quantum products for $Gr_2(\C^4)$.  For example, let us
calculate $a\circ a\circ a\circ a\circ a$.
From (5) we have $a\circ b\circ c = ae^t $, hence from
(2) and (3), $d\circ b=d\circ c= ae^t $.  It follows that
$a\circ a\circ a\circ a\circ a =2d\circ a\circ a$ (from [4])
$= 2d\circ a^2= 2d\circ (b+c)$ (from (1))  
$=4ae^t $ (from the previous calculation).  

In this example, the indirect approach just gives 
the subalgebra of the quantum cohomology algebra
of $Gr_2(\C^4)$ which is generated by the element $a\in H^2(Gr_2(\C^4);\C)$.
From the formulae above we have
$$
\align
a\circ a &=a^2\\
a\circ a\circ a &= a^3\\
a\circ a\circ a\circ a &= a^4 +2e^t \\
a\circ a\circ a\circ a\circ a &= 4ae^t 
\endalign
$$
and so we see that this algebra is generated abstractly by $a$ and $e^t $ with
the single relation $a^5-4ae^t =0$.  Note that this is neither the
quantum cohomology of  $Gr_2(\C^4)$ nor a subalgebra of the quantum cohomology of 
$\C P^5$. It is simply the algebra (i) in the Lefschetz decomposition, which
happens to be a subalgebra of the quantum cohomology of  $Gr_2(\C^4)$.

For other manifolds $M$ the situation is more complicated.  Even for
$Gr_2(\C^5)$, the algebra (i) fails to be closed under quantum multiplication.  
Another problem is that the induced homomorphism $\pi_2 M\to \pi_2 X$ may fail to be an
isomorphism. Nevertheless, it is possible to introduce a modified quantum product on $X$,
as explained in section 2 of \cite{Pa}, which reflects part of the quantum cohomology
of $M$. Using this,
the corresponding quantum differential equations and function $J$ can be defined.
If $M$ is a hypersurface (or complete intersection), similar behaviour to
that of the above example is expected.

In all known cases, as well as for the ambient manifold $X$ itself,
{\it the function $J$ turns out to be
related to a function of hypergeometric
\footnote{Recall that a power series $\sum a_n z^n$ is said to
be hypergeometric if $a_n/a_{n+1}$ is a rational function of $n$. 
We are obviously using the term rather loosely, as we have functions of
several variables whose values lie in the cohomology algebra.}
type,} in a specific way which
depends on the nature of the ambient space $X$ (see \cite{Gi5}
and \cite{BCKv} for the main examples).  That is, although the function $J$ was
originally expressed as a generating function for certain Gromov-Witten
invariants, it turns out to have an unexpectedly simple analytical form
--- as we have seen for $\C P^n$ and the examples in the appendices.

This surprising and
nontrivial fact, which reveals the hidden structure of the
Gromov-Witten \ll enumerative\rr data,  is often referred to as the Mirror
Theorem (or Mirror Identity). The name is a reference to the
Mirror Symmetry Conjecture, which implies that for each such $M$ there
exist a \ll mirror partner\rr  $M^\ast$ with the property that $J$ is related to
the integral of a holomorphic $n$-form over
a cycle in $M^\ast$ (a \ll period integral\rrr), which satisfies a Picard-Fuchs
differential equation.   The conjecture would predict (and to
some extent explain) the Mirror Theorem.

The model example where hypergeometric functions arise as period integrals
is the case of the two-dimensional real torus, a cubic curve in $\C P^2$.
The period integrals are classical hypergeometric functions.  In \cite{KLRY},
this case is examined as motivation, followed by (quartic) $K3$ surfaces in $\C P^3$, then 
(quintic) Calabi-Yau manifolds in $\C P^4$.  Although mirror pairs of
Calabi-Yau manifolds were the original focus, 
versions of Mirror Symetry for Fano manifolds were initiated
independently in \cite{EHX} and \cite{Gi4}. For a Fano
manifold $M$, however, the mirror partner is thought to be
(a partial compactification of) $(\C^\ast)^n$, where $n=\dim_{\C}M$,
rather than a manifold of the same kind as $M$.

A final brief comment on terminology may be helpful.
The term \ll Quantum Lefschetz Hyperplane Theorem\rr might be expected to
refer  to the relation between the quantum cohomology of a
hypersurface (or complete intersection) $M$ and the quantum cohomology
of an ambient projective space (or flag manifold or Fano toric manifold) $X$.
However, this term has so far been used (as in \cite{Ki2},\cite{BCKv})
to mean the Mirror Theorem for such $M$ and $X$.  This version of the
Mirror Theorem has been proved rigorously for various classes of $M$ and $X$,
and it justifies the enumerative predictions first made by physicists.  On the
other hand it does not by itself explain Mirror Symmetry; it could be
described as \ll Mirror Symmetry without the mirror manifolds\rrr.

\head
Appendix 1: Quantum differential equations for $F_{1,2}(\C^3)$
\endhead

{\eightpoint

We shall examine in detail the
quantum differential equations for the flag manifold
$$
F_3=F_{1,2}(\C^3) =
\{ (L,V) \in Gr_1(\C^3) \times Gr_2(\C^3) \st L\sub V \}.
$$
We use the notation of Appendix 1 of \cite{Gu}, in particular
the specific generators $a,b$ of the cohomology algebra.

Let us choose the ordered basis
$$
b_0=1,\ 
b_1=a,\ 
b_2=b,\ 
b_3=a^2,\ 
b_4=b^2,\ 
b_5=a^2b=ab^2=z
$$
of $H^\ast(F_3;\C)$, from which it follows that
the dual basis with respect to the intersection form is
$$
a_0=z,\ 
a_1=b^2,\ 
a_2=a^2,\ 
a_3=b,\ 
a_4=a,\ 
a_5=1.
$$
A general two-dimensional cohomology class will be denoted $t=t_1a+t_2b$.
We recall some further terminology from \cite{Gu}.   The (positive) K\"ahler
cone --- i.e. the cohomology classes representable by K\"ahler $2$-forms ---
consists of the classes $n_1 a + n_2 b$ with $n_1,n_2 > 0$. (The first Chern
class of the holomorphic tangent bundle is $2a+2b$.) The Mori cone ---
i.e. the homotopy classes (or their corresponding homology classes) which have 
holomorphic representatives --- are the classes
$D=d_1 B^2 + d_2 A^2$ with $d_1,d_2\ge 0$. All these facts may be found in
Appendix 1 of \cite{Gu}.

In \cite{Gu} the formal variables $q_1,q_2$ were defined initially by
$$
q^D = q^{d_1 B^2 + d_2 A^2} = 
(q^{B^2})^{d_1} (q^{A^2})^{d_2} = q_1^{d_1} q_2^{d_2}.   
$$
As we want to use the quantum product $\circ_t$ (rather than $\circ$)
we shall use the representation of $q_1,q_2$ as complex numbers given by
$q_1=e^{t_1}, q_2=e^{t_2}$, i.e.
$$
q^D = e^{\lan t,D \ran } =
e^{\lan t_1 b_1 + t_2 b_2, d_1 B^2 + d_2 A^2 \ran } = e^{t_1d_1 + t_2d_2}
= (e^{t_1})^{d_1} (e^{t_2})^{d_2}.
$$

The matrix-valued $1$-form $\om$ is 
$\om_t = M_1(t)dt_1+ M_2(t)dt_2$
where the matrices $M_1,M_2$ are (respectively) the
matrices of the quantum multiplication operators
$a\ct$, $b\ct$ on $H^\ast(F_3;\C)$. From the multiplication
rules in Appendix 1 of \cite{Gu} we have
$$
M_1(t) =
\pmatrix
0 & q_1 & 0 & 0 & 0 & q_1q_2 \\
1& 0 & 0 & 0 & 0 & 0 \\
0& 0 & 0 & q_1 & 0 & 0 \\
0& 1 & 1 & 0 & 0 & 0 \\
0& 0 & 1 & 0 & 0 & q_1 \\
0& 0 & 0 & 0 & 1 & 0 
\endpmatrix,
\quad
M_2(t) =
\pmatrix
0 & 0 & q_2 & 0 & 0 & q_1q_2\\
0& 0 & 0 & 0 & q_2 & 0 \\
1& 0 & 0 & 0 & 0 & 0 \\
0& 1 & 0 & 0 & 0 & q_2 \\
0& 1 & 1 & 0 & 0 & 0 \\
0& 0 & 0 & 1 & 0 & 0 
\endpmatrix.
$$

We are interested in the maps $H$ such that $\om = h dH H^{-1}$, i.e.
(invertible) matrix-valued functions whose column vectors
$\psi:H^2(F_3;\C)\to H^\ast(F_3;\C)$ form a basis of solutions of
the linear system
$$
h \frac{\b \psi}{\b t_1} = M_1 \psi,\quad 
h \frac{\b \psi}{\b t_2} = M_2 \psi.
$$
It turns out to be more convenient to work with the rows of $H$.
We write
$$
H=
\pmatrix
--- & J_0  & ---\\
--- & J_1 & ---\\
--- & J_2 & ---\\
--- & J_3 & ---\\
--- & J_4 & ---\\
--- & J_5 & ---
\endpmatrix,
$$
where $J_i:H^2(F_3;\C)\to H^\ast(F_3;\C)$. Let
$\bone=\frac{\b}{\b t_1}=q_1\frac{\b}{\b q_1}$, 
$\btwo=\frac{\b}{\b t_2} = q_2\frac{\b}{\b q_2}$.  Then the
system $h \bone H = M_1 H$ is equivalent to
$$
\align
h\bone J_0 &= q_1 J_1 + q_1q_2 J_5
\\
h\bone J_1 &= J_0
\\
h\bone J_2 &= q_1 J_3
\\
h\bone J_3 &= J_1 + J_2
\\
h\bone J_4 &= J_2 + q_1 J_5
\\
h\bone J_5 &= J_4
\endalign
$$
and the system $h \btwo H = M_2 H$ is equivalent to
$$
\align
h \btwo J_0 &= q_2 J_2 + q_1q_2 J_5
\\
h \btwo J_1 &= q_2 J_4
\\
h \btwo J_2 &= J_0
\\
h \btwo J_3 &= J_1 + q_2 J_5
\\
h \btwo J_4 &= J_1 + J_2
\\
h \btwo J_5 &= J_3.
\endalign
$$
Let $J=J_5$.   Then five of the above twelve equations may be used to express
$J_0,J_1,J_2,J_3,J_4$ in terms of $J$ as follows:
$$
\align
J_0 &= h^3 \bone^2 \btwo J - q_1 h \btwo J
\\
J_1 &= h^2 \bone\btwo J - h^2 \bone^2 J + q_1 J
\\
J_2 &= h^2 \bone^2 J - q_1 J
\\
J_3 &= h \btwo J
\\
J_4 &= h\bone J
\endalign
$$
(this may be done in various ways; we have made a particular choice).
The remaining seven equations reduce to the following system of equations
for $J$:
$$
\align
(h^2 \bone^2 + h^2\btwo^2 - h^2\bone\btwo -q_1 -q_2)J &=0 \tag 1
\\
(h^3 \bone \btwo^2 - h^3 \bone^2 \btwo - q_2 h \bone + q_1 h \btwo)J &=0 \tag 2
\\
(h^3 \bone^3 - q_1 h \bone - q_1 h \btwo )J&=q_1 h J \tag 3
\\
(h^4\bone^3 \btwo - 2q_1 h^2 \bone\btwo + 
q_1 h^2 \bone^2 - q_1^2 - q_1q_2 )J &=q_1 h^2 \btwo J \tag 4
\\
h^4 \bone^2\btwo^2 - q_1 h^2 \btwo^2 - q_2 h^2 \bone^2)J &=0. \tag 5
\endalign
$$
Equations $(3), (4), (5)$ follow from $(1), (2)$.    We conclude that the system
$h dH H^{-1} = \om$ is {\it equivalent} to the system
$$
D_1 J = 0,D_2 J = 0
$$
where
$$
H=
\pmatrix
--- & h^3 \bone^2\btwo J - q_1h \btwo J  & ---\\
--- & h^2\bone\btwo J - h^2\bone^2 J + q_1 J & ---\\
--- & h^2 \bone^2 J - q_1 J & ---\\
--- & h\btwo J & ---\\
--- & h\bone J & ---\\
--- & J & ---
\endpmatrix
$$
and
$$
\align
D_1 &= h^2 \bone^2 + h^2\btwo^2 - h^2\bone\btwo -q_1 - q_2\\
D_2 &= h^3 \bone \btwo^2 - h^3 \bone^2 \btwo - q_2 h\bone + q_1 h\btwo.
\endalign
$$
In the terminology
of \S 3, the equations $D_1J=0,D_2J=0$ are the quantum differential equations.

In general one can attempt to find formal power series solutions
to this kind of system.   In the case at hand, there is a particularly
simple \ll hypergeometric\rr formula, namely:
$$
J(t)=
\sum_{d_1,d_2\ge 0}
\frac
{\prod_{m=1}^{d_1+d_2} (a+b+mh)}
{\prod_{m=1}^{d_1} (a+mh)^3  \prod_{m=1}^{d_2} (b+mh)^3}
q_1^{a/h+d_1}
q_2^{b/h+d_2}
$$
(this is essentially the same as a function obtained in part III of \cite{Sc}
and in \cite{BCKv}).
The scalar (i.e. $H^0$-) component of this function is the hypergeometric function
$$
\sum_{d_1,d_2\ge 0}
\frac{1}{h^{2d_1+2d_2}}
\frac {(d_1+d_2)!} {(d_1!)^3 (d_2!)^3}
q_1^{d_1}q_2^{d_2}.
$$
Without addressing the question of
where this formula comes from, we just shall verify that $J$ is in fact a solution.

We have 
$h \bone q_1^{a/h+d_1} = (a+d_1h)q_1^{a/h+d_1}$,
$h \btwo q_2^{b/h+d_2} = (b+d_2h)q_2^{b/h+d_2}$.  Let us
write the above formula as 
$J(t)=\sum_{d_1,d_2\ge 0} a_{d_1,d_2} 
q_1^{a/h+d_1}
q_2^{b/h+d_2}$. Then the coefficient of $q_1^{a/h+d_1}
q_2^{b/h+d_2}$ in $D_1J$ is
$$
\spreadlines{3\jot}
\align
&a_{d_1,d_2} \{ (a+d_1h)^2 + (b+d_2h)^2 - (a+d_1h)(b+d_2h) \}
-a_{d_1-1,d_2} - a_{d_1,d_2-1}\\
&=a_{d_1,d_2} \left\{ (a+d_1h)^2 + (b+d_2h)^2 - (a+d_1h)(b+d_2h)
-
\frac
{(a+d_1h)^3+(b+d_2h)^3}
{a+b+(d_1+d_2)h}
\right\}
\\
&=\frac
{a_{d_1,d_2}}
{a+b+(d_1+d_2)h}
\{
(a+b+(d_1+d_2)h)(  (d_1^2+d_2^2-d_1d_2) h^2  +   \\
&\quad\quad\quad\quad\quad\quad  \quad\quad\quad\quad\quad\quad
 d_1(2a-b)h + d_2(2b-a) h )
-(a+d_1h)^3 - (b+d_2h)^3
\}
\endalign
$$
where we have used $a^2+b^2-ab=0$.
The coefficients of $h^3$ and $h^2$ (in the parentheses)
are both zero; the coefficient of $h$ is zero because of the
identity $a^2+b^2-ab=0$; the constant coefficient is zero
because of the identity $a^3+b^3=0$.

A similar verification works for $D_2$. The coefficient of 
$q_1^{a/h+d_1} q_2^{b/h+d_2}$ in $D_2J$ is
$$
\spreadlines{3\jot}
\align
&a_{d_1,d_2} \{ (a+d_1h)(b+d_2h)^2 - (a+d_1h)^2(b+d_2h) \}
-a_{d_1,d_2-1}(a+d_1h) + a_{d_1-1,d_2}(b+d_2h)\\
&=a_{d_1,d_2} \left\{ (a+d_1h)(b+d_2h)^2 - (a+d_1h)^2(b+d_2h)
-\frac
{(a+d_1h)(b+d_2h)^3-(a+d_1h)^3(b+d_2h)}
{a+b+(d_1+d_2)h}
\right\}
\\
&=\frac
{a_{d_1,d_2}(a+d_1h)(b+d_2h)}
{a+b+(d_1+d_2)h}
\{
(a+b+(d_1+d_2)h) ((b+d_2h) - (a+d_1h))
-(b+d_2h)^2+(a+d_1h)^2
\}
\endalign
$$
and the expression in the parentheses is zero.

According to the theorem at the end of \S 1 (see the beginning of \S 3
for the version involving $J$), 
the differential operators $D_1,D_2$ produce relations
$$
\gather
a\circ a + b\circ b - a\circ b - q_1 - q_2 = 0\\
a \circ b \circ b - a\circ a\circ b - q_2 a + q_1 b = 0
\endgather
$$
in the quantum cohomology algebra of $F_3$.   It turns out that these constitute
a complete set of relations (see Appendix 1 of \cite{Gu}).  So we
have found sufficiently many quantum differential equations 
in order to determine (up to isomorphism)
the quantum cohomology algebra.

The quantum product itself is not determined by these relations, but
it is of course determined by the function $H$ (because $H$ determines $\om$).   
The additional information involved is the form of the 
above expressions for the rows of $H$ in terms of $J$.  Our particular choices of
these expressions correspond to the quantum product formulae
$$
\align
a^2b &=
a\circ a\circ b - q_1 b \\
ab-a^2 &= a\circ b - a\circ a +q_1 \\
a^2 &= a\circ a - q_1 \\
b &= b\\
a &= a\\
1 &= 1.
\endalign
$$
In other words, we are taking
the elements of an additive basis of $H^\ast(F_3;\C)$ (written as cup products
of two-dimensional classes) and expressing them as quantum products of two-dimensional
classes.  It is easy to check directly
(using the computations in Appendix 1 of \cite{Gu}) 
that these particular quantum product formulae,
together with the two relations, determine the quantum products of
{\it arbitrary} cohomology classes.

To express this most conveniently, observe that we may factor the matrix $H$ as
$$
\align
H&=
\pmatrix
--- & h^3 \bone^2\btwo J - q_1h \btwo J  & ---\\
--- & h^2\bone\btwo J - h^2\bone^2 J + q_1 J & ---\\
--- & h^2 \bone^2 J - q_1 J & ---\\
--- & h\btwo J & ---\\
--- & h\bone J & ---\\
--- & J & ---
\endpmatrix\\
&=
\pmatrix
1 & 0 & 0 & -q_1 & 0 & 0 \\
0& 1 & 0 & 0 & 0 & q_1 \\
0& 0 & 1 & 0 & 0 & -q_1 \\
0& 0 & 0 & 1 & 0 & 0 \\
0& 0 & 0 & 0 & 1 & 0 \\
0& 0 & 0 & 0 & 0 & 1 
\endpmatrix
\pmatrix
--- & h^3 \bone^2\btwo J & ---\\
--- & h^2\bone\btwo J - h^2\bone^2 J & ---\\
--- & h^2 \bone^2 J & ---\\
--- & h\btwo J & ---\\
--- & h\bone J & ---\\
--- & J & ---
\endpmatrix
\quad
(= QH_0,\ \text{say}).
\endalign
$$
Then our observation is that the matrix $Q$ is the one which arises
when the basis elements $1, a, b, a^2, b^2, a^2b$ are expressed as
quantum polynomials in $a$ and $b$.

Finally we recall that the asymptotic version 
$H_{-\infty}: H^2(F_3;\C) \to \End(H^\ast(F_3;\C))$
(with the property $h d H_{-\infty} H_{-\infty}^{-1} = \om_{-\infty}$
where $(\om_{-\infty})_t(x)(y) = xy $) is given by
$H_{-\infty}(t):x\mapsto e^{t/h} x$.  The matrix of $H_{-\infty}(t)$
is therefore
$$
\pmatrix
{\ \ \ 1\ \ \ } & 0 & 0& 0& 0& 0 \\
\frac{t_1}h & {\ \ \ 1\ \ \ } & 0& 0& 0& 0 \\
\frac {t_2}h & 0& {\ \ \ 1\ \ \ }& 0& 0& 0 \\
\frac{2 t_1t_2 + t_1^2} {2h^2} & \frac{t_1+t_2} h & \frac {t_1} h &{\ \ \ 1\ \ \ } &0 &0 \\
\frac{2 t_1t_2 + t_2^2} {2h^2} & \frac{t_2} h & \frac{t_1+t_2} h &0 &{\ \ \ 1\ \ \ } &0 \\
\frac{t_1^2t_2 + t_1t_2^2} {2h^3} & \frac{2t_1t_2 + t_2^2} {2h^2} & 
\frac{2t_1t_2 + t_1^2} {2h^2} 
& \ \ \ \frac{t_2} h\ \ \  &\ \ \ \frac {t_1} h\ \ \  &{\ \ \ 1\ \ \ }
\endpmatrix
$$
The asymptotic version $J_{-\infty}$  of $J$ is the last row of this matrix.
One can verify that it satisfies the \ll classical\rr differential
equations
$$
\gather
( \bone^2 + \btwo^2 - \bone\btwo)J_{-\infty} = 0\\
( \bone \btwo^2 -  \bone^2 \btwo)J_{-\infty} = 0.
\endgather
$$
}

\head
Appendix 2: Quantum differential equations for $\Si_1$
\endhead

{\eightpoint

We shall compute explicit solutions of the
quantum differential equations for the Hirzebruch surface
$\Si_1 = \P(\Cal O(0)\oplus \Cal O(-1) )$, where $\Cal O(i)$
denotes the holomorphic line bundle on $\C P^1$ 
with first Chern class $i$. 
We use the notation of Appendix 2 of \cite{Gu}.

Let us choose the ordered basis
$$
b_0=1,\ 
b_1=x_1,\ 
b_2=x_4,\ 
b_3=x_1x_4 = z
$$
of the cohomology vector space. We then have
$$
a_0=z,\ 
a_1=x_2,\ 
a_2=x_1,\ 
a_3=1
$$
and $t=t_1x_1+t_2x_4$.

The  K\"ahler cone is spanned by $x_1,x_4$. (The first Chern class of
the tangent bundle is $2x_4+x_1$.)
The Mori cone (of holomorphically representable classes)
consists of all classes of the form $d X_1 + e X_2 $
with either $d\ge e\ge 0$ or $d=0, e\ge 0$.

In \cite{Gu} we defined the formal variables $q_1,q_2$
by 
$$
q^D = q^{d X_1 + e X_2} = 
(q^{X_1})^{d} (q^{X_2})^{e} = q_1^d q_2^e.
$$
Here we shall change notation to $r_1=q_2$, $r_2=q_1$, and (as in
Appendix 1) we shall consider these to be the complex numbers
$r_1=e^{t_1}$, $r_2=e^{t_2}$. Thus we have 
$$
q^D = e^{\lan t,D \ran } =
e^{\lan t_1 x_1 + t_2 x_4,d X_1 + e X_2 \ran} = e^{t_1e + t_2d}
= (e^{t_1})^{e} (e^{t_2})^{d} = r_1^e r_2^d.
$$

The matrix-valued $1$-form $\om$ is 
$\om_t = M_1(t)dt_1+ M_2(t)dt_2$
where the matrices $M_1,M_2$ are (respectively) the
matrices of the quantum multiplication operators
$x_1\ct$, $x_4\ct$ on $H^\ast(\Si_1;\C)$. From the multiplication
rules in Appendix 2 of \cite{Gu} we have
$$
M_1(t) =
\pmatrix
0 & 0 & 0 & r_1r_2 \\
1 & -r_1 & 0 & 0 \\
0 & r_1 & 0 & 0 \\
0 & 0 & 1 & 0
\endpmatrix,
\quad
M_2(t) =
\pmatrix
0 & 0 & r_2 & r_1r_2 \\
0 & 0 & 0 & r_2\\
1 & 0 & 0 & 0 \\
0 & 1 & 1 & 0
\endpmatrix.
$$

We are interested in the maps $H$ such that $\om = h dH H^{-1}$. As in 
Appendix 1 we shall rewrite this equation using the notation
$$
H=
\pmatrix
--- & J_0  & ---\\
--- & J_1 & ---\\
--- & J_2 & ---\\
--- & J_3 & ---
\endpmatrix,
$$
where $J_i:H^2(F_3;\C)\to H^\ast(F_3;\C)$.  Let
$\bone=\frac{\b}{\b t_1}=r_1\frac{\b}{\b r_1}$, 
$\btwo=\frac{\b}{\b t_2} = r_2\frac{\b}{\b r_2}$. Then the
system $h \bone H = M_1 H$ is equivalent to
$$
\align
h\bone J_0 &= r_1r_2 J_3
\\
h\bone J_1 &= J_0 - r_1 J_1
\\
h\bone J_2 &= r_1 J_1
\\
h\bone J_3 &= J_2
\endalign
$$
and the system $h \btwo H = M_2 H$ is equivalent to
$$
\align
h \btwo J_0 &= r_2 J_2 + r_1r_2 J_3
\\
h \btwo J_1 &= r_2 J_3
\\
h \btwo J_2 &= J_0
\\
h \btwo J_3 &= J_1 + J_2.
\endalign
$$
Let $J=J_3$.   Then three of the above eight equations may be used to express
$J_0,J_1,J_2$ in terms of $J$ as follows:
$$
\align
J_0 &= h^2 \bone\btwo  J 
\\
J_1 &= h\btwo J - h\bone J
\\
J_2 &= h \bone J
\endalign
$$
The remaining five equations reduce to the following system of equations
for $J$:
$$
\align
(h^2 \btwo^2  - h^2 \bone\btwo  - r_2) J &= 0 \tag 1 \\
(h^2 \bone^2  - r_1 h \btwo  + r_1 h \bone) J&=0 \tag 2\\
(h^3\bone^2\btwo  - r_1r_2) J &=0 \tag 3\\
(h^3 \bone\btwo^2  - r_2 h \bone  - r_1r_2) J &=0 \tag4
\endalign
$$
Equations $(3), (4)$ follow from $(1), (2)$.  Hence the original system
$hdHH^{-1} = \om$ is {\it equivalent} to the system
$$
D_1 J = 0,D_2 J = 0
$$
where
$$
H=
\pmatrix
--- & h^2\bone \btwo J & --- \\
--- & h\btwo J- h\bone J & --- \\
--- & h\bone J & ---\\
--- & J & ---
\endpmatrix
$$
and
$$
\align
D_1 &= h^2 \btwo^2 - h^2\bone\btwo -  r_2\\
D_2 &= h^2 \bone^2 - r_1 h \btwo + r_1 h \bone.
\endalign
$$

An explicit  solution of the quantum differential equations 
$D_1 J = 0$, $D_2 J = 0$ is given by
$$
J(t)=
\sum_{d,e\ge 0}
\frac
{1}
{\prod_{m=1}^{e} (x_1+mh)^2 \prod_{m=1}^{d} (x_4+mh)}
\frac
{\prod_{m=-\infty}^{0} (x_2+mh)}
{\prod_{m=-\infty}^{d-e} (x_2+mh)}
r_1^{x_1/h+e}
r_2^{x_4/h+d}
$$
(this function was obtained in \cite{Gi5}; further information can
be found in chapter 11 of \cite{Co-Ka}).
We shall verify this fact directly.
We have 
$h \bone r_1^{x_1/h+e} = (x_1+eh) r_1^{x_1/h+e}$,
$h \btwo r_2^{x_4/h+d} = (x_4+dh) r_2^{x_4/h+d}$.  Let us
write $J(t)=\sum_{e,d\ge 0} a_{e,d} 
r_1^{x_1/h+e}
r_2^{x_4/h+d}$. Then the coefficient of $r_1^{x_1/h+e}
r_2^{x_4/h+d}$ in $D_1J$ is
$$
\align
&a_{e,d} \{ (x_4+dh)^2 - (x_1+eh)(x_4+dh) \} - a_{e,d-1}\\
&=a_{e,d} \{ (x_4+dh)^2 - (x_1+eh)(x_4+dh) 
- (x_4+dh)(x_2 + (d-e)h )\}\\
&= a_{e,d}(x_4+dh)\{  x_4+dh - (x_1+eh) - (x_2 + (d-e)h) \}\\
&=  a_{e,d}(x_4+dh)(  x_4 - x_1 -x_2) \\
&= 0
\endalign
$$
since $x_4=x_1+x_2$.

Similarly, the coefficient of $r_1^{x_1/h+e}
r_2^{x_4/h+d}$ in $D_2J$ is
$$
\align
&a_{e,d} (x_1+eh)^2 - (x_4+dh)a_{e-1,d} + (x_1+eh)a_{e-1,d}\\
&=a_{e,d} \{ (x_1+eh)^2  -  
\frac{(x_4+dh)(x_1+eh)^2}{x_2+(d-e)h}  + 
\frac{(x_1+eh)^3}{x_2+(d-e)h}  
\}\\
&= a_{e,d} (x_1+eh)^2 \{
(x_2+(d-e)h) -(x_4+dh) + (x_1+eh) \}\\
&= a_{e,d} (x_1+eh)^2 
(  x_2 - x_4 +x_1) \\
&= 0
\endalign
$$
again, as $x_4=x_1+x_2$.

As in Appendix 1,
the differential operators $D_1,D_2$ produce two relations
$$
\gather
x_4\circ x_4 - x_1\circ x_4 -r_2 =0\\
x_1\circ x_1 - r_1 x_4 + r_1 x_1 = 0
\endgather
$$
in the quantum cohomology algebra of $\Si_1$.   Again it turns out that these constitute
a complete set of relations (see Appendix 2 of \cite{Gu}). 

The quantum product is determined by these relations together with
the specific form of $H$.  One may verify
that the latter is equivalent to the formulae
$$
\align
x_1x_4 &= x_1\circ x_4 \\
x_4-x_1 &= x_4-x_1 \\
x_1&= x_1\\
1 &= 1,
\endalign
$$
i.e. to the single formula $x_1\circ x_4=x_1x_4$.
From the computations in Appendix 2 of \cite{Gu}, one sees
that this particular quantum product formula,
together with the two relations, determine the quantum products of
{\it arbitrary} cohomology classes.  

In contrast to the case of $F_3$, the \ll quantization matrix\rr $Q$
is the identity here.

The asymptotic version of $H$, i.e. the operator $e^{t/h}$, is in this
case
$$
H_{-\infty}(t)=
\pmatrix
1 & 0 & 0 & 0 \\
\frac{t_1}h & 1 & 0 & 0 \\
\frac{t_2}h & 0 & 1 & 0 \\
\frac{t_2^2 + 2t_1t_2}{2h^2} & \frac{t_2}h & \frac{t_1+t_2}h & 1
\endpmatrix.
$$
The asymptotic version  of $J$, i.e. the last row of this matrix,
satisfies the \ll classical\rr differential equations
$$
\gather
(\btwo^2 - \bone\btwo)J_{-\infty}=0\\
\bone^2 J_{-\infty} =0.
\endgather
$$

}

\eightpoint \it

\no  Department of Mathematics
\newline
Graduate School of Science
\newline
Tokyo Metropolitan University
\newline
Minami-Ohsawa 1-1, Hachioji-shi
\newline
Tokyo 192-0397, Japan

\no martin\@math.metro-u.ac.jp

\enddocument